\documentclass[12pt]{article}

\usepackage{amsthm}
\usepackage{amsmath}
\usepackage{mathrsfs}
\usepackage{stmaryrd}
\usepackage{amssymb}
\usepackage{diagbox}
\usepackage{tikz-cd}
\usepackage{amsfonts}
\usepackage{mathtools}
\usepackage{authblk}
\usepackage{bm}
\usepackage[utf8]{inputenc}
\usepackage{indentfirst}
\usepackage{enumerate}
\usepackage[top=2.54cm, bottom=2.54cm, left=3.28cm, right=3.28cm]{
geometry}
\usepackage{verbatim}
\usepackage[colorlinks,linkcolor=blue]{hyperref}

\newcommand{\inv}{^{-1}}

\newcommand{\tm}{\tilde{M}}

\DeclareMathOperator{\sign}{sign}

\newtheorem{thm}{Theorem}
\theoremstyle{definition}

\newtheorem{lem}{Lemma}

\newtheorem{cor}{Corollary}
\newtheorem{exm}{Example}

\makeatletter
\newcommand{\proofpart}[2]{
  \par
  \addvspace{\medskipamount}
  \noindent\textbf{Part #1. #2}
  \par\nobreak\smallskip
  \@afterheading
}
\makeatother

\begin{document}

\title{The Generation of All Regular Rational Orthogonal Matrices}
\author{\small Quanyu Tang$^{{\rm a}}$\quad\quad Wei Wang$^{\rm a}$\thanks{Corresponding author: wang\_weiw@xjtu.edu.cn}\quad\quad Hao Zhang$^{\rm b}$
\\
{\footnotesize$^{\rm a}$School of Mathematics and Statistics, Xi'an Jiaotong University, Xi'an 710049, P. R. China}\\
{\footnotesize$^{\rm b}$School of Mathematics, Hunan University, Changsha 410082, P. R. China}}

\date{}
	\maketitle
	\begin{abstract}
A \emph{rational orthogonal matrix} $Q$ is an orthogonal matrix with rational entries, and $Q$ is called \emph{regular} if each of its row sum equals one, i.e., $Qe = e$ where $e$ is the all-one vector. This paper presents a method for generating all regular rational orthogonal matrices using the classic Cayley transformation. Specifically, we demonstrate that for any regular rational orthogonal matrix $Q$, there exists a permutation matrix $P$ such that $QP$ does not possess an eigenvalue of $-1$. Consequently, $Q$ can be expressed in the form $Q = (I_n + S)^{-1}(I_n - S)P$, where $I_n$ is the identity matrix of order $n$, $S$ is a rational skew-symmetric matrix satisfying $Se = 0$, and $P$ is a permutation matrix. Central to our approach is a pivotal intermediate result, which holds independent interest: given a square matrix $M$, then $MP$ has $-1$ as an eigenvalue for every permutation matrix $P$ if and only if either every row sum of $M$ is $-1$ or every column sum of $M$ is $-1$.

	\end{abstract}

\section{Introduction}

\subsection{Why regular rational orthogonal matrices?}
An $n \times n$ matrix $Q$ is an \emph{orthogonal matrix} if $Q^{\rm T}Q = I_n$, where $I_n$ is the identity matrix of order $n$. $Q$ is called a \emph{rational orthogonal matrix} if it is an orthogonal matrix whose entries are rational, and it is \emph{regular} if the sum of the entries in each of its rows (and columns) equals one, i.e., $Qe = e$ where $e = (1,1,\dots,1)^{\rm T}$ represents the all-one vector.

Orthogonal matrices are fundamental in various fields of mathematics. This paper focuses on rational orthogonal matrices, which are notably prevalent in graph theory and combinatorics. Recall that a \emph{Hadamard matrix} $H=(h_{ij})$ is a matrix with $h_{ij}\in{\{1,-1\}}$ such that $H^{\rm T}H=nI_n$. If $H$ is a Hadamard matrix of order $n=m^2$, then $\frac{1}{m}H$ is a rational orthogonal matrix. Similarly, a \emph{conference matrix} $C = (c_{ij})$, characterized by $c_{ii} = 0$ and $c_{ij} \in \{1, -1\}$ for $i \neq j$ such that $C^{\rm T}C = (n-1)I_n$, becomes a rational orthogonal matrix as $\frac{1}{m}C$ for $n = m^2 + 1$.

\begin{exm} The following matrices $H$ and $C$ are regular rational orthogonal matrices of orders $4$ and $10$, respectively.
$$H=\frac{1}{2}\left[ \begin{array}{cccc}
-1&1&1&1\\
1&-1&1&1\\
1&1&-1&1\\
1&1&1&-1
\end{array}\right], $$

$$C=\frac{1}{3}\left(
\begin{array}{cccccccccc}
 0 & 1 & 1 & 1 & 1 & 1 & 1 & -1 & -1 & -1 \\
 1 & 0 & 1 & 1 & 1 & -1 & -1 & 1 & 1 & -1 \\
 1 & 1 & 0 & -1 & -1 & 1 & 1 & 1 & 1 & -1 \\
 1 & 1 & -1 & 0 & 1 & -1 & 1 & 1 & -1 & 1 \\
 1 & 1 & -1 & 1 & 0 & 1 & -1 & -1 & 1 & 1 \\
 1 & -1 & 1 & -1 & 1 & 0 & 1 & -1 & 1 & 1 \\
 1 & -1 & 1 & 1 & -1 & 1 & 0 & 1 & -1 & 1 \\
 -1 & 1 & 1 & 1 & -1 & -1 & 1 & 0 & 1 & 1 \\
 -1 & 1 & 1 & -1 & 1 & 1 & -1 & 1 & 0 & 1 \\
 -1 & -1 & -1 & 1 & 1 & 1 & 1 & 1 & 1 & 0 \\
\end{array}
\right).$$

\end{exm}

Our interest in regular rational orthogonal matrices is primarily driven by their applications in spectral graph theory, a field that employs matrices and linear algebra tools to address problems in graph theory and its applications (see, for example,~\cite{BH,CDS}). A graph $G$ is said to be \emph{determined by its spectrum} (abbreviated as DS) if any graph sharing the same spectrum as $G$ is isomorphic to it. The question of ``Which graphs are DS?" remains a significant unresolved query in spectral graph theory (refer to~\cite{DH1,DH2}). Demonstrating that a graph is DS often presents substantial challenges. In their research, Wang and Xu~\cite{WX1} suggest approaching the spectral characterization problem through the lens of generalized spectrum.

Let $G$ be a graph with an adjacency matrix $A$ on $n$ vertices. A graph $G$ is \emph{determined by its generalized spectrum} (abbreviated as DGS) if, for any graph $H$, whenever $G$ and $H$ are cospectral with cospectral complements, then $H$ must be isomorphic to $G$. Define the \emph{walk-matrix} of graph $G$ as $W(G):=[e, Ae, \ldots, A^{n-1}e]$. In the works of Wang~\cite{Wang1, Wang2}, the following theorem is established:

\begin{thm}[Wang~\cite{Wang1,Wang2}]\label{Wang}
If $2^{-\lfloor n/2 \rfloor}{\det W(G)}$ (which is always an integer) is odd and square-free, then $G$ is DGS.
 \end{thm}

The main idea behind the proof of Theorem~\ref{Wang} is articulated in the subsequent theorem, which establishes a connection between the generalized cospectrality of graphs and regular rational orthogonal matrices.

\begin{thm}[Johnson and Newman~\cite{JN}; Wang and Xu~\cite{WX1}]\label{rational}  Two graphs $G$ and $H$ are cospectral with cospectral complements if and only if there exists an orthogonal matrix $Q$ such that
$Q^{\rm T}A(G)Q=A(H).$
Moreover, if $W(G)$ is nonsingular, then the above $Q$ is a regular rational orthogonal matrix which is unique and equals $W(G)W(H)^{-1}$.
\end{thm}

According to Theorem~\ref{rational}, a graph $G$ with $\det W(G) \neq 0$ is determined by its generalized spectrum (DGS) if and only if for all regular rational orthogonal matrices $Q$, the condition $Q^{\rm T}AQ$ yielding a (0,1)-matrix implies that $Q$ must be a permutation matrix.

Regular rational orthogonal matrices have also been instrumental in the construction of cospectral graphs, which are graphs that share identical spectra. Notably, the GM-switching method, devised by Godsil and McKay~\cite{GM}, has been demonstrated to be an effective approach for constructing cospectral graphs  with cospectral complements, which uses a regular rational orthogonal matrices of a particular type. An analogue of the GM-switching was given in~\cite{QJW}, which uses regular orthogonal matrices of another type.

Therefore, the study of regular rational orthogonal matrices not only advances our understanding of generalized spectral characterizations of graphs but also facilitates the construction of cospectral graphs.

\subsection{Cayley transformation and rational orthogonal matrices}

Let $S$ be a skew-symmetric matrix, i.e., $S^{\rm T} = -S$. It can be readily confirmed that the matrix
\begin{equation}\label{EQ1}
Q = (I_n + S)^{-1}(I_n - S)
\end{equation}
is orthogonal. The mapping $\Phi$ from the set of all skew-symmetric matrices to the set of all orthogonal matrices, defined by the transformation $$S \longmapsto Q = (I_n + S)^{-1}(I_n - S),$$ is known as the \emph{Cayley transformation}. It is clear that the matrix $Q$ as defined in Eq.~\eqref{EQ1} does not possess an eigenvalue of $-1$. Conversely, if $-1$ is not an eigenvalue of $Q$, then the matrix $S = (I_n + Q)^{-1}(I_n - Q)$ is skew-symmetric. Therefore, the Cayley transformation $\Phi$ is injective but not surjective. Liebeck and Osborn~\cite{LO} devised a method to generate all rational orthogonal matrices using the Cayley transformation. They showed, among others, that for any rational orthogonal matrix $Q$, there always exists a diagonal matrix $D$ with all diagonal entries $\pm1$ such that $DQ$ does not have an eigenvalue $-1$.

\begin{thm} [Liebeck and Osborn~\cite{LO}]\label{thm1} Let $Q$ be an orthogonal matrix. Then there exists a diagonal matrix $D$ with all diagonal entries $\pm 1$ such that $-1$ is not an eigenvalue of $DQ$. Moreover, every rational orthogonal matrix $Q$ can be expressed as $Q=D(I_n+S)^{-1}(I_n-S)$ for some rational skew-symmetric matrix $S$ and some
diagonal orthogonal matrix $D$.

\end{thm}

\subsection{Main result}
The aforementioned theorem by Liebeck and Osborn does not extend to regular rational orthogonal matrices, primarily because the matrix $DQ$ may no longer retain its regularity. It is noted that the matrix $Q$, as defined in Eq.~\eqref{EQ1}, is regular if and only if the vector sum of each row of $S$ is zero, i.e., $Se = 0$. Inspired by a similar concept in~\cite{LO}, this paper develops a method to generate all regular rational orthogonal matrices via Cayley transformation. Instead of employing a diagonal orthogonal matrix $D$, we utilize a permutation matrix $P$ and ensure that every row sum of $S$ is zero. We demonstrate that for any regular rational orthogonal matrix $Q$, there is always a permutation matrix $P$ such that $QP$ does not possess an eigenvalue of $-1$.

The main result of our study is encapsulated in the following theorem:

\begin{thm}\label{main} Let $Q$ be any $n\times n$ regular rational orthogonal matrix. Then there exists a permutation matrix $P$ such that $QP$ does not have an eigenvalue $-1$. Moreover, every regular rational orthogonal matrix $Q$ can be expressed as $Q=(I_n+S)^{-1}(I_n-S)P$ for some rational skew-symmetric matrix $S$ with $Se=0$, and some permutation matrix $P$.
\end{thm}

Let $RO_n(\mathbb{Q})$ denote the set of all regular orthogonal matrices of order $n$ over $\mathbb{Q}$, and $SS_n(\mathbb{Q})$ the set of all skew-symmetric matrices of the same order over $\mathbb{Q}$. Let $\Sigma_n$ represent the set of all permutation matrices of order $n$. Then we have

\begin{cor}Using the notations above, we have
 $$RO_n(\mathbb{Q})=\{(I_n+S)^{-1}(I_n-S)P:\, P\in{\Sigma_n}, S\in{SS_n(\mathbb{Q})}~{\rm with}~ Se=0\}.$$
\end{cor}

The crux of the proof of Theorem~\ref{main} hinges on another theorem, which is noteworthy on its own.

\begin{thm}\label{thm:mainthm}
Let $M$ be an $n\times n$ matrix. If $\det(I_n+MP)= 0$ for all permutation matrix $P$, then either $Me= -e$ or $M^{\rm T}e= -e$.
\end{thm}
Assuming at the moment that Theorem~\ref{thm:mainthm} is true, we give the proof of Theorem~\ref{main}:
\begin{proof}[Proof of Theorem~\ref{main}] Note that $Qe=e$. It follows that $Qe\neq -e$ and $Q^{\rm T}e\neq -e$. According to Theorem~\ref{thm:mainthm}, there exists a
permutation matrix $P$ such that $QP$ does not have an eigenvalues $-1$. Let $S=(I+QP)^{-1}(I-QP)$. Then $S$ is a skew-symmetric matrix with rational entries
and $Se=0$. Consequently, $Q=(I+S)^{-1}(I-S)P^{\rm T}\in{RO_n(\mathbb{Q})}$.

\end{proof}

The rest of the paper is devoted to the proof of Theorem~\ref{thm:mainthm}.

\section{Proof of Theorem~\ref{thm:mainthm}}
We begin by establishing some notations. Let $n > 1$ be a positive integer. Define $[n] = \{1, 2, \dots, n\}$ and let $S_n$ denote the symmetric group on $[n]$. Consider $M = (m_{ij})$ as an $n \times n$ matrix with column vectors $M_1, \dots, M_n$ and row vectors $\tm_1^{\rm T}, \dots, \tm_n^{\rm T}$ respectively. For any subsets $I = \{i_1, \dots, i_k\}$ and $J = \{j_1, \dots, j_k\}$ of $[n]$, where $i_1 < \dots < i_k$ and $j_1 < \dots < j_k$, we define $M_{I,J} = (m_{ij})_{i \in I, j \in J}$ as the submatrix of $M$ formed by row indices $I$ and column indices $J$. The determinant of this submatrix is denoted by $|M_{I,J}|$. Additionally, for any $\sigma \in S_n$, we use ${\rm sign}(\sigma, I)$ to denote the sign of $\sigma$ restricted to $I$, corresponding to the permutation $\begin{pmatrix} i_1 & i_2 & \cdots & i_k \\ \sigma(i_1) & \sigma(i_2) & \cdots & \sigma(i_k) \end{pmatrix}$. For any $\sigma \in S_n$, let $P_{\sigma}$ be the permutation matrix with nonzero entries at positions $(\sigma(i), i)$, for $i = 1, 2, \dots, n$.

\begin{lem}\label{lem:sum0}
    Let $I,J\subseteq [n]$ with $|I|=|J|\geq 2$, then
    \[\sum_{\sigma\in S_n,\sigma(I)=J}\sign(\sigma,I)=0.\]
\end{lem}

\begin{proof}
    Fix $j_1,j_2\in J$, and put $\tau=(j_1\ j_2)\in S_n$. Then $J=\sigma(I)$ if and only if $J=\tau\sigma(I)$. Moreover, we have $\sign(\sigma,I)=-\sign(\tau\sigma,I)$. So
    \[2\sum_{\sigma\in S_n,\sigma(I)=J}\sign(\sigma,I)=\sum_{\sigma\in S_n,\sigma(I)=J}\sign(\sigma,I)+\sum_{\sigma\in S_n,\tau\sigma(I)=J}\sign(\tau\sigma,I)=0.\]
\end{proof}

\begin{lem}\label{lem:enMen}
Let $M$ be an $n\times n$ matrix. If $|I_n+MP|= 0$ for all permutation matrix $P$, then $e^{\rm T}Me=-n$, i.e., the sum of all entries of $M$ is $-n$.
\end{lem}

\begin{proof}
It is well-known that
\[|\lambda I_n+M|=\sum_{k=0}^n\left(\sum_{I\subseteq [n]\atop |I|=k}|M_{I,I}|\right)\lambda^{n-k},\]
where $M_{\emptyset,\emptyset}=1$ by convention. In particular, we have
\[|I_n+M|=\sum_{I\subseteq [n]}|M_{I,I}|.\]
Now for any $\sigma\in S_n$ and $I=\{i_1,\dots,i_k\}\subseteq [n]$ with $i_1<\dots<i_k$, recall that $\sign(\sigma,I)=1$ if the inversion number of $(\sigma(i_1),\dots,\sigma(i_k))$ is even, and $\sign(\sigma,I)=-1$ if the inversion number of $(\sigma(i_1),\dots,\sigma(i_k))$ is odd. Let $P_\sigma$ be the permutation matrix corresponding to $\sigma$, and $M=(M_1,M_2,\dots,M_n)$, then
\[MP_{\sigma}=(M_{\sigma(1)},M_{\sigma(2)},\dots,M_{\sigma(n)}).\]
So for any $I\subseteq [n]$, we have $|(MP_{\sigma})_{I,I}|=\sign(\sigma,I)|M_{I,\sigma(I)}|$.
So we get
\begin{equation}\label{eq:imp}
\begin{split}
    0=&\sum_{\sigma\in S_n}|I_n+MP_{\sigma}|=\sum_{\sigma\in S_n}\sum_{I\subseteq [n]}|(MP_{\sigma})_{I,I}|=\sum_{\sigma\in S_n}\sum_{I\subseteq [n]}\sign(\sigma,I)|M_{I,\sigma(I)}|\\
    =&\sum_{I\subseteq [n]}\sum_{J\subseteq [n]\atop |I|=|J|}\sum_{\sigma\in S_n\atop J=\sigma(I)}\sign(\sigma,I)|M_{I,J}|.
\end{split}
\end{equation}
For fixed $I,J\subseteq [n]$, if $|I|=|J|=0$ or $1$, then $\sign(\sigma,I)=1$ for all $\sigma\in S_n$, and hence
\[\sum_{\sigma\in S_n\atop J=\sigma(I)}\sign(\sigma,I)=(n-|I|)!.\]
When $|I|=|J|\geq 2$, by the Lemma \ref{lem:sum0}, we have
\[\sum_{\sigma\in S_n, J=\sigma(I)}\sign(\sigma,I)=0.\]
So from Eq.~\eqref{eq:imp}, we get $\sum_{i,j} m_{ij}=-n$.
\end{proof}

\begin{lem}\label{lem:mij}
    Let $M,M_i$ and $\tm_j$ be defined as before. Suppose that $n+\sum_{i,j=1}^nm_{ij}=0$. Then for any $1\leq i,j\leq n$, we have
    \begin{equation}\label{eq:sigmaij}
        \sum_{\sigma\in S_n\atop \sigma(i)=j}|I_n+MP_{\sigma}|=-(n-2)!(1+e^{\rm T}\tm_i)(1+e^{\rm T}M_j).
    \end{equation}
\end{lem}

\begin{proof}
We first consider the case $i=j=n$. Then the same calculation as in Lemma \ref{lem:enMen} gives
\[\sum_{\sigma\in S_n\atop \sigma(n)=n}|I_n+MP_{\sigma}|=\sum_{I\subseteq [n]}\sum_{J\subseteq [n]\atop |I|=|J|}\left(\sum_{\substack{\sigma\in S_n\\ \sigma(I)=J\\ \sigma(n)=n}}\sign(\sigma,I)\right)|M_{I,J}|:=\sum_{I\subseteq [n]}\sum_{J\subseteq [n]\atop |I|=|J|}C(I,J)|M_{I,J}|.\]
If $n\in I$ and $n\notin J$, then there is no contribution of $\sigma$ in $C(I,J)$, i.e., $C(I,J)=0$. Similarly, if $n\notin I$ and $n\in J$, we also have $C(I,J)=0$. So we only need to consider the cases $n\in I\cap J$ and $n\notin I\cup J$. We first suppose that $n\in I\cap J$.

If $|I|=|J|=1$, i.e., $I=J=\{n\}$, then $\sign(\sigma,I)=1$ for all $\sigma\in S_n$ with $\sigma(n)=n$. So we have $C(I,J)=(n-1)!$.

If $|I|=|J|=2$, then $I=\{i,n\}$ and $J=\{j,n\}$ for some $1\leq i,j\leq n-1$. In this case, for any $\sigma\in S_n$ with $\sigma(i)=j$ and $\sigma(n)=n$, we always have $\sign(\sigma,I)=1$, so we have $C(I,J)=(n-2)!$. Moreover, in this case we have
\[|M_{I,J}|=\begin{vmatrix}
    m_{ij}&m_{in}\\m_{nj}&m_{nn}
\end{vmatrix}=m_{ij}m_{nn}-m_{in}m_{nj}.\]

If $|I|=|J|\geq 3$, then by Lemma~\ref{lem:sum0}, we see that
\[C(I,J)=\sum_{\sigma\in S_{n-1}\atop \sigma(I\backslash\{n\})=J\backslash\{n\}}\sign(\sigma,I\backslash\{n\})=0.\]

Next we suppose that $n\notin I\cup J$. Similar to the discussion as above, when $|I|=|J|=0$, we have $C(\emptyset,\emptyset)=(n-1)!$. When $|I|=|J|=1$, say $I=\{i\},J=\{j\}$, then we have $C(I,J)=(n-2)!$ and $|M_{I,J}|=m_{ij}$. By Lemma \ref{lem:sum0} again, we have $C(I,J)=0$ when $|I|=|J|\geq 2$.

Combining all the information together, we get
\begin{equation}\label{eq:inmp}
\begin{split}
    \sum_{\sigma\in S_n\atop \sigma(n)=n}|I_n+MP_{\sigma}|=&(n-1)!m_{nn}+(n-2)!\sum_{1\leq i,j\leq n-1}(m_{ij}m_{nn}-m_{in}m_{nj})\\
    &+(n-1)!+(n-2)!\sum_{1\leq i,j\leq n-1}m_{ij}.
\end{split}
\end{equation}
To simplify the notation, we put $x=\sum_{i=1}^{n-1}m_{in}$ and $y=\sum_{j=1}^{n-1}m_{nj}$. By our assumption, we have $\sum_{1\leq i,j\leq n}m_{ij}=-n$, so $\sum_{1\leq i,j\leq n-1}m_{ij}=-n-m_{nn}-x-y$. So from Eq. \eqref{eq:inmp}, we get
\begin{align*}
    &\sum_{\sigma\in S_n, \sigma(n)=n}|I_n+MP_{\sigma}|\\
    =&(n-2)!\left((n-1)m_{nn}+(1+m_{nn})(-n-m_{nn}-x-y)-xy+n-1\right)\\
    =&-(n-2)!(1+m_{nn}+x)(1+m_{nn}+y).
\end{align*}
Now we consider the general case, for any $1\leq i,j\leq n$, we assume that $i\neq n$. Fix a permutation $\tau$ such that $\tau(n)=j$ and $\tau(i)=n$. Put $M'=P_{\tau}MP_{\tau}$, then the sum of entries of $M'$ is again $-n$, so by the proof above, we have
\begin{align*}
    -(n-2)!(1+e^{\rm T}\tm_n')(1+e^{\rm T}M_n')=&\sum_{\sigma\in S_n, \sigma(n)=n}|I_n+M'P_\sigma|\\
    =&\sum_{\sigma\in S_n, \sigma(n)=n}|P_{\tau}||I_n+MP_\tau P_\sigma P_\tau||P_{\tau\inv}|\\
    =&\sum_{\sigma\in S_n, \sigma(n)=n}|I_n+M P_{\tau\sigma\tau} |.
\end{align*}
Notice that when $\sigma$ runs through $S_n$ which fixes $n$, $\tau\sigma\tau$ runs through all the permutations in $S_n$ which sends $i$ to $j$. So the last summation is precisely the LHS of Eq.~\eqref{eq:sigmaij}. On the other hand, we note that the $n$-th row of $P_\tau MP_\tau$ is the $\tau\inv(n)$-th row of $MP_{\tau}$, so $e^{\rm T}\tm'_n=e^{\rm T}M_i$. Similarly, we have $e^{\rm T}M_n'=e^{\rm T}M_j$.

Finally, if $i=n$, since we have shown that the Eq.~\eqref{eq:sigmaij} holds for $i=j=1$, we can fix a permutation $\tau$ such that $\tau(1)=j$ and $\tau(n)=1$, then by repeating the proof above, we can get the desired result.
\end{proof}

\begin{proof}[Proof of Theorem \ref{thm:mainthm}]
Suppose that $Me\neq -e$ and $M^{\rm T}e\neq e$, then there exist two indices $i$ and $j$ such that $e^{\rm T}\tm_i\neq -1$ and $e^{\rm T}M_j\neq -1$. Then by Lemma~\ref{lem:mij}, we have
\[0=\sum_{\sigma\in S_n\atop \sigma(i)=j}|I_n+MP_{\sigma}|=-(n-2)!(1+e^{\rm T}\tm_i)(1+e^{\rm T}M_j)\neq 0,\]
which gives a contradiction.
\end{proof}

\section*{Acknowledgments}
The research of the second author is supported by National Key Research and Development Program of China 2023YFA1010203, National Natural Science Foundation of China (Grant No.\,12371357) and the third author is supported by Fundamental Research Funds for the Central Universities (Grant No.\,531118010622).


\begin{thebibliography}{22}

		
	\bibitem{BH} A.E. Brouwer and W.H. Haemers, Spectra of Graphs, \emph{Springer}, 2012.

\bibitem{CDS}	D. M. Cvetkovi\'{c}, M. Doob, H. Sachs, Spectra of Graphs, \emph{Academic Press,
	NewYork}, 1982.

\bibitem{DH1}
	E. R. van Dam, W. H. Haemers, Which graphs are determined by their spectrum?
	\emph{Linear Algebra Appl}. 373 (2003) 241-272.
	
	\bibitem{DH2}
	E. R. van Dam, W. H. Haemers, Developments on spectral characterizations
	of graphs, \emph{Discrete Math.}, 309 (2009) 576-586.	
	
	
	\bibitem{JN}C. R. Johnson, M. Newman, A note on cospectral graphs, \emph{J. Combin. Theory, Ser. B}, 28 (1980) 96-103.

\bibitem{GM} C.D. Godsil, B.D. McKay, Constructing cospectral graphs, \emph{Aequ. Math.},
	25 (1982) 257-268.
	

\bibitem{LO} H. Liebeck, A. Osborn, The Generation of All Rational Orthogonal Matrices, \emph{Amer. Math. Monthly}, Vol. 98, No. 2. (Feb., 1991), pp. 131-133.	
	
	

\bibitem{QJW} L. Qiu, Y. Ji, W. Wang, On a theorem of Godsil and McKay on the construction of cospectral graphs, \emph{Linear Algebra Appl}. 603 (2020) 265-274.
	
	
	
	
	\bibitem{WX1}	W. Wang, C. X. Xu, A sufficient condition for a family of graphs being determined
	by their generalized spectra, \emph{Eur. J. Combin.} 27 (2006) 826-840.
	
		
	
	\bibitem{Wang1} W. Wang, Generalized spectral characterization revisited, \emph{Electron. J.
	Combin.} 20 (4) (2013), $\sharp$ P4.
	
	
	\bibitem{Wang2}  W. Wang, A simple arithmetic criterion for graphs being determined by their generalized
spectra, \emph{J. Combin. Theory, Ser. B}, 122 (2017) 438-451.
	
	
	
	
	
\end{thebibliography}
\end{document}